\documentclass{amsart}
\usepackage{setspace}

\usepackage{amsmath}
\usepackage{amsthm}
\usepackage{amssymb}

\usepackage[latin1]{inputenc} 

\newtheorem{theorem}{Theorem}[section]

\newtheorem{proposition}[theorem]{Proposition}
\newtheorem{corollary}[theorem]{Corollary}
\newtheorem{remark}[theorem]{Remark}

\newtheorem{definition}[theorem]{Definition}

\numberwithin{equation}{section}

\begin{document}
\title[Properties of a type of the entropy of an ideal and the divergence ]
{Some properties of a type of the entropy of an ideal and the divergence of two ideals}

\author{Nicu\c{s}or Minculete}
\address{Faculty of Mathematics and Computer Science, Transilvania University\\
 Iuliu Maniu street 50, Bra\c{s}ov 500091, Romania}
\email{minculete.nicusor@unitbv.ro}

\author{Diana Savin}
\address{Faculty of Mathematics and Computer Science, Transilvania University\\
 Iuliu Maniu street 50, Bra\c{s}ov 500091, Romania}
\email{diana.savin@unitbv.ro; dianet72@yahoo.com}

\date{}
\begin{abstract}
The aim of this paper is to study certain properties of the Kullback--Leibler distance between two positive integer numbers or between two ideals. We present some results related the entropy of a positive integer number and the divergence of two numbers. We also study the entropy of some types of ideals and the divergence of two ideals. Finally, we find some inequalities, involving the  entropy  $H$ of an  exponential divisor of a positive integer, respectively the  entropy  $H$ of an  exponential divisor of an ideal.
\end{abstract}
\maketitle
\section{Introduction and Preliminaries} 
\noindent 
\smallskip\\

The entropy is defined, in information theory, as a measure of uncertainty. One of the most used entropies is the Shannon entropy ($H_S$), which is given for a probability distribution ${\bf p}=\{p_1,...,p_r\}$ thus
$$H_S({\bf p})= - \sum_{i=1}^{r}p_i\cdot \log p_i,$$
where $\sum_{i=1}^{r}p_i=1$.
An useful property of the Shannon entropy is the additivity,  $H_S({\bf pq})=H_S({\bf p})+H_S({\bf q})$, where ${\bf p}=\{p_1,...,p_r\}$, ${\bf q}=\{q_1,...,q_r\}$ and ${\bf pq}=\{p_1q_1,...,p_1q_r,...,p_rq_1,...,p_rq_r\}$,
and the recursivity,
$$
H_S(p_1,p_2,...,p_r)=H_S(p_1+p_2,p_3,...,p_r)+(p_1+p_2)H_S(\frac{p_1}{p_1+p_2},\frac{p_2}{p_1+p_2}).
$$

In \cite{CT} the relative entropy or (Kullback-Leibler distance) between two probability distributions ${\bf p}=\{p_1,...,p_r\}$  and ${\bf q}=\{q_1,...,q_r\}$  was introduced in the following way:
$$
D({\bf p}||{\bf q}):=-\sum_{i=1}^{r}p_i\cdot \log \frac{q_i}{p_i}=\sum_{i=1}^{r}p_i\cdot \log \frac{p_i}{q_i},
$$ 
where $\sum_{i=1}^{r}p_i=\sum_{i=1}^{r}q_i=1$.

Let $n\in\mathbb{N},$  $n\geq2,$ according to the Fundamental Theorem of Arithmetic, $n$ is written uniquely
 $n=p^{\alpha_{1}}_{1}p^{\alpha_{2}}_{2}...p^{\alpha_{r}}_{r}$, where $r,\alpha_{1}, \alpha_{2},...,\alpha_{r}\in\mathbb{N}^{*}$ and  $p_{1}, p_{2},...,p_{r}$ are distinct prime positive integers. We denote by $\Omega\left(n\right)=\alpha_{1}+\alpha_{2}+...+\alpha_{r}$ and $p\left(\alpha_{i}\right)=\frac{\alpha_{i}}{\Omega\left(n\right)},$
$\left(\forall\right)$ $i=\overline{1,r}$. In \cite{Minculete1}, Minculete and  Pozna introduced the notion of entropy of  $n$  like this: 

\begin{equation}
H\left(n\right)= - \sum_{i=1}^{r}p\left(\alpha_{i}\right)\cdot \log \: p\left(\alpha_{i}\right),    \tag{1.1}
\end{equation}
where $\log$ is the natural logarithm. We have $H(1)=0$ (by convention).

In \cite{Minculete1}, the authors found an equivalent form for the entropy of $n,$ namely:
\begin{equation}
H\left(n\right)=\log \: \Omega\left(n\right) - \frac{1}{\Omega\left(n\right)}\cdot \sum_{i=1}^{r}\alpha_{i}\cdot \log \: \alpha_{i}.  \tag{1.2}
\end{equation}
Also, they proved the following result:
\begin{proposition}
\label{onedotone}
\begin{equation}
0\leq H\left(n\right)\leq \log \: \omega\left(n\right), \; \left(\forall\right) \; n\in \mathbb{N},\; n\geq 2,  \tag{1.3}
\end{equation}
where $ \omega\left(n\right)$  is the number of distinct prime factors of $n.$
\end{proposition}
We are giving some examples.\\
If $n=10=2\cdot 5,$ we have:
$$ H\left(10\right)=\log \: 2 -  \frac{1}{2}\cdot 2\cdot \log \: 1=\log \: 2=0.6931.... $$
If $n=100=2^{2}\cdot 5^{2},$ we have:
$$ H\left(100\right)=\log \: 4 -  \frac{1}{4}\cdot 2\cdot 2\cdot \log \: 2=\log \: 2=0.6931.... $$
If $n=2^{3},$ we have:
$$ H\left(8\right)=\log \: 3 -  \frac{1}{3}\cdot 3\cdot \log \: 3=0.$$
If if $n=40=2^{3}\cdot 5,$ we have:
$$ H\left(40\right)=\log \: 4 -  \frac{1}{4}\cdot 3\cdot \log \: 3= \frac{1}{4}\cdot \log\left(\frac{4^4}{3^3}\right)=2.2493.... $$
In  \cite{Minculete1}, the authors proved that:

\begin{remark}
\label{onedottwo}
i) If $n =p^{\alpha}$, with $\alpha$ a positive integer and $p$ a prime positive integer, then $H\left(n\right)=0;$\\
ii) If $n =p_{1}\cdot p_{2}\cdot...\cdot p_{r},$ with $p_{1}, p_{2},..., p_{r}$ distinct prime positive integers, then $H\left(n\right)=log\:\omega\left(n\right);$\\
iii) If $n =\left(p_{1}\cdot p_{2}\cdot...\cdot p_{r}\right)^{\alpha},$ with $\alpha$ a positive integer and $p_{1}, p_{2},..., p_{r}$ distinct prime positive integers, then $H\left(n\right)=log\:\omega\left(n\right).$
\end{remark}
We remark that $H(n^{\alpha})=H(n)$, for any positive integer $\alpha$.\\

From \cite{MS} we see that if $\gcd(m,n)=1$, then $$H(mn)\neq H(m)+H(n).$$ Here is studied the relationship between $H(mn)$ and $H(m)+H(n)$, where $m,n\in\mathbb{N}^*$, $m,n\geq 2$, $n=p^{\alpha_{1}}_{1}p^{\alpha_{2}}_{2}...p^{\alpha_{r}}_{r}$ and $m=q^{\beta_{1}}_{1}q^{\beta_{2}}_{2}...q^{\beta_{s}}_{s}$, $r,s\in\mathbb{N}^*$.

Among the important results obtained in \cite{MS} we mention some below. If $\gcd(n,p)=1$, with $p$ is a prime number and $n,\alpha\in\mathbb{N}^*$, then we have
It follows that 
\begin{equation}\label{ms_5}
H(np^\alpha)=\frac{\Omega(n)H(n)}{\Omega(n)+\alpha}+\log(\Omega(n)+\alpha)-\frac{\Omega(n)\log \Omega(n)+\alpha\log\alpha}{\Omega(n)+\alpha}.
\end{equation}
Using the definition of $H$, we deduce the equality
\begin{equation}\label{eq}
H(mn)-H(m)-H(n)=\frac{\Omega(n)}{\Omega(m)(\Omega(m)+\Omega(n))}\sum_{i=1}^r\alpha_i\log \alpha_i
\end{equation}
$$+\frac{\Omega(m)}{\Omega(n)(\Omega(m)+\Omega(n))}\sum_{j=1}^s\beta_j\log \beta_j-\log\frac{\Omega(m)\Omega(n)}{\Omega(m)+\Omega(n)}.
$$
From relation \eqref{eq} and the definition of $H$, we have 
\begin{align}\label{eq.1}
& H(mn)-H(m)-H(n)\\
&=\frac{\Omega(n)\log(\Omega(m))+\Omega(m)\log(\Omega(n))}{\Omega(n)+\Omega(m)}-\frac{\Omega(n)H(m)+\Omega(m)H(n)}{\Omega(n)+\Omega(m)}-\log\frac{\Omega(m)\Omega(n)}{\Omega(m)+\Omega(n)}.\nonumber
\end{align}

We assume that $m=p^kq$ and $n=p^kt$, where $p,q,t$ are distinct prime numbers and $k\in\mathbb{N}^*$. Then the inequality 
$$H(mn)<H(m)+H(n)$$
holds.
If $m=p_1^kp_2$ and $n=q_1^kq_2$, where $p_1,p_2,q_1,q_2$ are distinct prime numbers and $k\in\mathbb{N}^*$, then we deduce the inequality 
$$H(mn)\geq H(m)+H(n).$$
Equality holds for $k=1$.

In \cite{Minculete1} it is studied the Kullback--Leibler distance between two positive integer numbers. The entropy concept is used in several developments like the distance between two random variables or like the relative information between two distributions.  

We consider two positive integer numbers $n=p^{\alpha_{1}}_{1}p^{\alpha_{2}}_{2}...p^{\alpha_{r}}_{r}$ and $m=q^{\beta_{1}}_{1}q^{\beta_{2}}_{2}...q^{\beta_{r}}_{r}$. Next, we define the probabilities $p(\alpha_i)=\frac{\alpha_i}{\Omega(n)}$  and  $p(\beta_i)=\frac{\beta_i}{\Omega(m)}$, for every $i\in\{1,2,...,r\}$, with $\sum_{i=1}^{r}p(\alpha_i)=1$ and $\sum_{i=1}^{r}p(\beta_i)=1$.

Similar with the Kullback--Leibler distance between two probability distributions, was defined the Kullback--Leibler distance between two positive integer numbers.

In \cite{Minculete1}, the Kullback--Leibler distance between two positive integer number $n$ and $m$ ($n,m\geq 2$) with $\omega(n)=\omega(m)$  and factorization thus  $n=p^{\alpha_{1}}_{1}p^{\alpha_{2}}_{2}...p^{\alpha_{r}}_{r}$ and $m=q^{\beta_{1}}_{1}q^{\beta_{2}}_{2}...q^{\beta_{r}}_{r}$ is given as                                             
$$
D(n||m):=-\sum_{i=1}^{r}p(\alpha_i)\cdot \log \frac{p(\beta_i)}{p(\alpha_i)}.
$$
This can be written in following way
$$
D(n||m)=\log\frac{\Omega(m)}{\Omega(n)}-\frac{1}{\Omega(n)}\sum_{i=1}^{r}\alpha_i\cdot \log \frac{\beta_i}{\alpha_i}.
$$
For example, if $n=100=2^25^2$  and $m=200=2^35^2$, then the  divergence of $n$ and $m$ is $D(n||m)=\log\frac{5}{4}-\frac{1}{2}\log\frac{3}{2}=\frac{1}{2}\log\frac{25}{24}$=0.0088....

It is easy to see that if $n=m$, then we have $D(n||m)=0$. 
If we have the integer numbers $n$ and $m$ with $\omega(n)=\omega(m)$  and factorization thus  $n=p^{\alpha_{1}}_{1}p^{\alpha_{2}}_{2}...p^{\alpha_{r}}_{r}$ and $m=q^{\alpha_{1}}_{1}q^{\alpha_{2}}_{2}...q^{\alpha_{r}}_{r}$, then we obtain $D(n||m)=D(m||n)=0$. Generally, we have  $D(n||m)\neq D(m||n)$.

Subbarao  introduced the notion of \textit{exponential divisor} of a positive integer and he found some properties of these divisors (see \cite{Straus}, \cite{Subbarao}). So, if $n$ is a positive integer,  $n>1,$ it can be written uniquely
as $n=p_1^{\alpha_1}p_2^{\alpha_2}\cdot\cdot\cdot p_r^{\alpha_r}$, where $r\in\mathbb{N}^{*},$ $p_{1}, p_{2},...,p_{r}$ are distinct prime positive integers and $\alpha_{1}, \alpha_{2},...,\alpha_{r}\in\mathbb{N}^{*}.$
The positive integer $d=p_1^{\beta_1}p_2^{\beta_2}\cdot\cdot\cdot p_r^{\beta_r}$ (with $\beta_{1}, \beta_{2},...,\beta_{r}\in\mathbb{N}^{*}$) is called \textit{exponential divisor} or \textit{e-divisor} of $n=p_1^{\alpha_1}p_2^{\alpha_2}\cdot\cdot\cdot p_r^{\alpha_r}$, 
if $\beta_i|\alpha_i$, for every $i\in \{1,...,r\}.$
The number of exponential divisors of $n$ is denoted by $\tau^{\left( e\right)}\left(n\right)$ and if $n>1$ denoted as above,  we have $\tau^{\left( e\right)}\left(n\right)=\tau\left( \alpha_1\right)\tau\left( \alpha_2\right)\cdot\cdot\cdot\tau\left( \alpha_r\right),$
where $\tau\left( \alpha_i\right)$ is the number of natural divisors of $ \alpha_i,$  $\left(\forall\right)$ $i=\overline{1,r}.$ By convention $\tau^{\left( e\right)}\left(1\right)=1.$ We note $d|_{(e)}n$ .\\
\indent Many properties of the number of the exponential divisors of a positive integer $n$ can be found in the articles \cite{Minculete}, \cite{Min_Sav}, \cite{Pet}, \cite{Sandor_2}, \cite{Sandor}, \cite{Toth}, \cite{Wu}. When $d|_{(e)}n$, we have  $\omega(n)=\omega(d)$, so, we deduce                                           
\begin{equation}\label{div_1}
D(n||d)=\log\frac{\Omega(d)}{\Omega(n)}-\frac{1}{\Omega(n)}\sum_{i=1}^{r}\alpha_i\cdot \log \frac{\beta_i}{\alpha_i}.
\end{equation}
Since $-\frac{1}{\Omega(n)}\sum_{i=1}^{r}\alpha_i\cdot \log \frac{\beta_i}{\alpha_i}=\frac{1}{\Omega(n)}\sum_{i=1}^{r}\alpha_i\cdot \log \frac{\alpha_i}{\beta_i}\geq 0$, we deduce
$$D(n||d)\geq\log\frac{\Omega(d)}{\Omega(n)}.$$

For $\gamma(n)=p_{1}p_{2}...p_{r}$, which is the lowest divisor of $n$, we calculate the divergence between $n$ and $\gamma(n)$ and we obtain
\begin{equation}\label{div_2}
D(n||\gamma(n))=\log\frac{\omega(n)}{\Omega(n)}+\frac{1}{\Omega(n)}\sum_{i=1}^{r}\alpha_i\log \alpha_i.
\end{equation}
Taking into account the definition of entropy we find the following relation:
\begin{equation}\label{div_2.1}
D(n||\gamma(n))+H(n)=\log\omega(n).
\end{equation}
From \eqref{div_1} and \eqref{div_2}, we deduce the difference
$$
D(n||d)-D(n||\gamma(n))=\log\frac{\Omega(d)}{\omega(n)}-\frac{1}{\Omega(n)}\sum_{i=1}^{r}\alpha_i\log \beta_i.
$$
 
In this paper we study some properties of the Kullback--Leibler distance between two positive integer numbers or between two ideals. In Section 2 we present some results related the entropy of a positive integer number and the divergence of two numbers. In Section 3 we study the entropy of some types of ideals and the divergence of two ideals. In the last section, we find some inequalities, involving the  entropy  $H$ of an  exponential divisor of a positive integer, respectively the  entropy  $H$ of an  exponential divisor of an ideal.
  
\section{Some results related the entropy of a positive integer number and the divergence of two numbers}

\noindent 
\smallskip\\
\begin{proposition}
\label{prop2.1}
For two natural numbers $m,n\geq 2$ with $\omega(m)=\omega(n)$, we have
\begin{equation}\label{div_2.2}
D(n||m)=H(m)-H(n)+\sum_{i=1}^{r}\left(\frac{\beta_i}{\Omega(m)}-\frac{\alpha_i}{\Omega(n)}\right)\log \beta_i.
\end{equation}
\end{proposition}
\begin{proof}
From the definition of divergence, we find
\begin{align*}
D(n||m)&=\log\frac{\Omega(m)}{\Omega(n)}-\frac{1}{\Omega(n)}\sum_{i=1}^{r}\alpha_i\cdot \log \frac{\beta_i}{\alpha_i}\\
&=\log\Omega(m)-\frac{1}{\Omega(n)}\sum_{i=1}^{r}\alpha_i \log \beta_i-H(n)\nonumber\\
&=H(m)+\frac{1}{\Omega(m)}\sum_{i=1}^{r}\beta_i \log \beta_i-\frac{1}{\Omega(n)}\sum_{i=1}^{r}\alpha_i \log \beta_i-H(n)\nonumber
\end{align*} 
Therefore, we prove the relation of the statement.
\end{proof}
\begin{remark}\label{remark2.2}
From relation \eqref{div_2.2}, we deduce the equalities:
\begin{equation}\label{div-en}
D(n||\gamma(n))=H(\gamma(n))-H(n)=\log\omega(n)-H(n)
\end{equation}
and
\begin{align*}
D(np^{\alpha}||np^{\beta})&=H(np^{\beta})-H(np^{\alpha})\\
&+\frac{\left(\alpha-\beta\right)\Omega(n)\left(\log\Omega(n)-H(n)\right)}{\left(\Omega(n)+\alpha\right)\left(\Omega(n)+\beta\right)}+\frac{\beta\log\beta}{\Omega(n)+\beta}-\frac{\alpha\log\beta}{\Omega(n)+\alpha}\\
&=H(np^{\beta})-H(np^{\alpha})\\
&+\frac{\left(\alpha-\beta\right)\Omega(n)\left(\log\Omega(n)-H(n)-\log\beta\right)}{\left(\Omega(n)+\alpha\right)\left(\Omega(n)+\beta\right)}.
\end{align*}
\end{remark}
\begin{proposition}
\label{prop2.2}
For $n\geq 3$, we have
\begin{equation}\label{div_2.3}
D(n||\gamma(n))+H(n)\leq\log\log n-\log\log\log n+\log c_1,
\end{equation}
where $c_1=1.38402...$.
\end{proposition}
\begin{proof}
In \cite{Ro}, G. Robin proved that $\omega(n)\leq\frac{\log n}{\log\log n}c_1$, for all $n\geq 3$ and $c_1=1.38402...$. Using this relation and relation \eqref{div_2.1}, we obtain the desired relation. 
\end{proof}
\begin{proposition}\label{ms.d}
If $\gcd(n,p)=1$, with $p$ a prime number and $n,\alpha\in\mathbb{N}^*$, then we have
\begin{equation}\label{sec2_1d}
\lim_{p\to\infty}\lim_{\alpha\to\infty}[\log\left(\omega(n)+1\right)-D(np^{\alpha}||\gamma(n)p)]=0. 
\end{equation}
\end{proposition}
\begin{proof}
From relation \eqref{div-en} we deduce $D(np^{\alpha}||\gamma(n)p)=\log\omega(np^{\alpha})-H(np^{\alpha})=\log\left(\omega(n)+1\right)-H(np^{\alpha})$ and from \cite[Proposition 3]{MS} we have $\lim_{p\to\infty}\lim_{\alpha\to\infty}H(np^\alpha)=0$. Therefore, we prove that the limit of the statement is true.
\end{proof}
\begin{proposition}\label{ms.2.d}
We assume that $m=p^kq$ and $n=p^kt$, where $p,q,t$ are distinct prime numbers and $k\in\mathbb{N}^*$. Then the inequality 
$$D(m||\gamma(m))+D(n||\gamma(n))<\log\frac{4}{3}+D(mn||\gamma(mn))$$
holds.
\end{proposition}
\begin{proof}
We apply relation \eqref{div-en} twice and taking into account that $H(mn)<H(m)+H(n)$, from \cite[Proposition 4]{MS}, we find the statement.
\end{proof}
\begin{proposition}
\label{ms.3.d}
We assume that $m=p_1^kp_2$ and $n=q_1^kq_2$, where $p_1,p_2,q_1,q_2$ are distinct prime numbers and $k\in\mathbb{N}^*$. Then we have the following inequality 
$$D(m||\gamma(m))+D(n||\gamma(n))\geq D(mn||\gamma(mn)).$$
Equality holds for $k=1$.
\end{proposition}
\begin{proof}
For $k=1$, we deduce that $m=p_1p_2$ and $n=q_1q_2$, which implies $D(m||\gamma(m))=D(n||\gamma(n))=D(mn||\gamma(mn))=0$.
For $k\geq 2$, we apply relation \eqref{div-en} twice and taking into account that $H(mn)\geq H(m)+H(n)$, from \cite[Proposition 5]{MS}, we find the inequality of the statement.
\end{proof}
\begin{theorem}
\label{ms.4.1.d}
Let $m,n$ be two natural numbers such that $\gcd(m,n)=1$ and $D(m||\gamma(m))\leq \log \frac{\omega(m)}{2}$, $D(n||\gamma(n))\leq \log \frac{\omega(n)}{2}$. Then the following inequality 
$$D(m||\gamma(m))+D(n||\gamma(n))\leq D(mn||\gamma(m)\gamma(n)).$$
holds.
\end{theorem}
\begin{proof}
Since $\gcd(m,n)=1$, we deduce that $\omega(mn)=\omega(m)+\omega(n)$ and $\gamma(mn)=\gamma(m)\gamma(n)$. Again, we apply relation \eqref{div-en} twice and taking into account that $H(mn)\leq H(m)+H(n)$, from \cite[Theorem 1]{MS}, we find the inequality of the statement.
\end{proof}
	
If $n$ is a number $k$-free (see \cite{Sandor_2}), so 
$n=p^{\alpha_{1}}_{1}p^{\alpha_{2}}_{2}...p^{\alpha_{r}}_{r}>1$, with $\alpha_i\leq k-1$, where $k\geq 2$, then using the entropy of $n$, we find (see \cite{Minculete1})
$$
\log\Omega(n)-\frac{\omega(n)}{\Omega(n)}\left(k-1\right)\log\left(k-1\right)\leq H(n)\leq \log\omega(n).
$$
From relation \eqref{div-en} and from above relation we prove
$$
0\leq D(n||\gamma(n))\leq\log\frac{\omega(n)}{\Omega(n)}+\frac{\omega(n)}{\Omega(n)}\left(k-1\right)\log\left(k-1\right).
$$
\section{The entropy of some types of ideals and the divergence of two ideals}

\noindent 
\smallskip\\

Let $K$ be an algebraic number field and let $\mathcal{O}_{K}$ be the ring of integers of $K$. The set of prime ideals of the ring $\mathcal{O}_{K}$ is denoted by Spec$\left(\mathcal{O}_{K}\right).$ Let $I$ be an ideal of the ring $\mathcal{O}_{K}.$ According to the Fundamental Theorem of Dedekind rings, there exist and they are unique $g\in\mathbb{N}^{*},$ the distinct ideals $P_{1},$ $P_{2},$..., $P_{g}\in Spec\left(\mathcal{O}_{K}\right)$ and the distincy numbers $e_{1},$ $e_{2},$...,$e_{g}\in\mathbb{N}^{*}$ such that $I=P^{e_{1}}_{1}\cdot P^{e_{2}}_{2}\cdot...\cdot P^{e_{g}}_{g}.$ 
In the paper \cite{MS}, we introduced the notion of entropy of an ideal of a ring of algebraic integers $\mathcal{O}_{K}$ in the following way:
\begin{definition}
\label{twodotone}
(Definition 1 from \cite{MS}).  Let $I\neq \left(0\right)$ be an ideal of the ring  $\mathcal{O}_{K}$, decomposed as above. We define the entropy of the ideal $I$ as follows:
\begin{equation}\label{sec3_3}
H\left(I\right):= - \sum_{i=1}^{g}\frac{e_{i}}{\Omega(I)} \log \: \frac{e_{i}}{\Omega(I)},     
\end{equation}
where $ \Omega\left(I\right) =  e_{1}+ e_{2}+...+ e_{g}.$
\end{definition}

In the article \cite{MS} we obtained the following equivalent form, for the entropy of the ideal $I:$
\begin{equation}\label{sec3_4}
H\left(I\right)=\log \: \Omega\left(I\right) - \frac{1}{\Omega\left(I\right)}\cdot \sum_{i=1}^{g}e_{i}\cdot \log \: e_{i}.  
\end{equation}
Also, in \cite{MS} we studied the entropy of ideals of the form $p\mathcal{O}_{K},$ with $p$ a prime integer.\\
\begin{definition}
\label{twodotone.1}
Let $I,J\neq \left(0\right)$ be two ideals of the ring  $\mathcal{O}_{K}$, decomposed as $I=P^{e_{1}}_{1}\cdot P^{e_{2}}_{2}\cdot...\cdot P^{e_{g}}_{g}$ and $J=Q^{f_{1}}_{1}\cdot Q^{f_{2}}_{2}\cdot...\cdot Q^{f_{g}}_{g}$. We define the divergence of the ideals $I$ and $J$ as follows:
\begin{equation}\label{sec3_3.1}
D\left(I||J\right):= \log\frac{\Omega(J)}{\Omega(I)}-\frac{1}{\Omega(I)}\sum_{i=1}^{g}e_i\cdot \log \frac{f_i}{e_i},     
\end{equation}
where $ \Omega\left(I\right) =  e_{1}+ e_{2}+...+ e_{g}$ and $ \Omega\left(J\right) =  f_{1}+ f_{2}+...+ f_{g}.$
\end{definition}
For the ideal $I=P^{e_{1}}_{1}\cdot P^{e_{2}}_{2}\cdot...\cdot P^{e_{g}}_{g}$, we denoted by $\gamma(I)=P_1P_2...P_g$ and $\omega(I)=g$, similar to numbers. Therefore we obtain
$$
D\left(I||\gamma(I)\right)= \log\frac{\omega(I)}{\Omega(I)}-\frac{1}{\Omega(I)}\sum_{i=1}^{g}e_i\cdot \log \frac{1}{e_i}=\log\omega(I)-H\left(I\right).     
$$

We denote by $n$ the degree of the extension of fields $\mathbb{Q}\subset K,$ where  $n\in\mathbb{N}$, $n\geq2.$
According to the fundamental theorem of Dedekind rings, 
the ideal $p\mathcal{O}_{K} $ is written uniquely in the form:
$$p\mathcal{O}_{K}=P^{e_{1}}_{1}\cdot P^{e_{2}}_{2}\cdot...\cdot P^{e_{g}}_{g},$$
where  $g\in\mathbb{N}^*,$   $e_{1}, e_{2},..., e_{g}\in\mathbb{N}^*$  and $P_{1},$ $P_{2},$..., $P_{g}\in Spec\left(\mathcal{O}_{K}\right).$
 The number $e_{i}$ ($i = \overline{1,g}$) is called the ramification index of $p$ at the ideal $P_{i}$. 
The following result is known (see \cite{ireland},  \cite{Ribenboim}):

\begin{proposition}
\label{onedotthree}
In the above notations, we have:\\
i)
$$\sum_{i=1}^{g}e_{i}f_{i}=[K :\mathbb{Q}] = n,$$
where  $f_{i}$ is the residual degree of $p,$ meaning $f_{i}=\left[\mathcal{O}_{K}/P_{i} : \mathbb{Z}/p\mathbb{Z}\right],$ $i = \overline{1,g}.$\\ 
ii) If moreover $\mathbb{Q} \subset K$ is a Galois extension, then $e_{1}=e_{2}=...=e_{g}$ (denoted by $e$), $f_{1}=f_{2}=...=f_{g}$ (denoted by $f$).
Therefore, $efg=n.$
\end{proposition}
In \cite{MS} we found the following properties of the entropy of ideals of the form $p\mathcal{O}_{K},$ with $p$ a prime integer.
\begin{remark}
\label{twodotfour}
Let $K$ be an algebraic number field and  let $\mathcal{O}_{K}$ be its  ring of integers.  Let $p$ be a prime positive integer. If $p$  is inert or totally ramified in the ring  $\mathcal{O}_{K},$  then $H\left(p\mathcal{O}_{K}\right)=0.$

\end{remark}

\begin{proposition}
\label{twodotfive}
 Let $n$ be a positive integer, $n\geq2$ and let $p$ be a positive prime integer. Let $K$ be an algebraic number field of degree $[K :
\mathbb{Q}] = n$ and let $\mathcal{O}_{K}$ be its the ring of integers. Then:\\
 \begin{equation}\label{sec3_5}
0\leq H\left(p\mathcal{O}_{K}\right)\leq \log \: \omega\left(p\mathcal{O}_{K}\right)\leq  \log \:n, \; 
\end{equation}
where $ \omega\left(p\mathcal{O}_{K}\right)$  is the number of distinct prime factors of the ideal $p\mathcal{O}_{K}.$
\end{proposition}

\begin{proposition}
\label{twodotsix}
Let $K$ be an algebraic number field and let $\mathcal{O}_{K}$ be its the ring of integers. Let $p$ be a prime positive integer. If the extension of fields $\mathbb{Q}\subset K$ is a Galois extension, then
$$H\left(p\mathcal{O}_{K}\right)=\log \: \omega\left(p\mathcal{O}_{K}\right).$$
\end{proposition}

We are studying what happens to the divergence of two ideals of the ring of algebraic integers of an algebraic nunber field $K$ when the extension of fields $\mathbb{Q}\subset K$ is a Galois extension a Galois extension. 

\begin{proposition}
\label{twodotseven}
Let $K$ be an algebraic number field such that  $\mathbb{Q}\subset K$ is a Galois extension. Let $\mathcal{O}_{K}$ be its ring of algebraic  integers and let $p$  and $q$  be two distinct prime positive integers such that  $\omega\left(p\mathcal{O}_{K}\right)= \omega\left(q\mathcal{O}_{K}\right).$ Then,
the divergence 
$$D\left(p\mathcal{O}_{K}||q\mathcal{O}_{K}\right)=0.$$
\begin{proof}
We denote by $g=\omega\left(p\mathcal{O}_{K}\right)= \omega\left(q\mathcal{O}_{K}\right).$ Since the extension of fields $\mathbb{Q}\subset K$ is a Galois extension, applying Proposition \ref{onedotthree}.,  it results that 
$p\mathcal{O}_{K}=P^{e}_{1}\cdot P^{e}_{2}\cdot...\cdot P^{e}_{g}$ and $J=Q^{e^{'}}_{1}\cdot Q^{e^{'}}_{2}\cdot...\cdot Q^{e^{'}}_{g}.$ Then $ \Omega\left(p\mathcal{O}_{K}\right) =  ge$ and $ \Omega\left(q\mathcal{O}_{K}\right) =  ge^{'}.$
According to Definition \ref{twodotone.1}, we have:
$$D\left(p\mathcal{O}_{K}||q\mathcal{O}_{K}\right):= \log\frac{e^{'}}{e}-\frac{1}{eg}\cdot  eg\cdot \log \frac{e^{'}}{e}=0.$$

\end{proof}

\end{proposition}
We are dealing with the entropy of the general ideals of a ring of algebraic integers $\mathcal{O}_{K}.$\\
We give the following example: let the quadratic field $K=\mathbb{Q}\left(i \sqrt{19}\right).$ Its class number is $1$ and its ring of  of algebraic integers $\mathcal{O}_{K}= \mathbb{Z}\left[\frac{1+i\sqrt{19}}{2}\right].$
It results that this is a principal ring, so it is a Dedekind ring. Let the ideal
 $$I=35\mathbb{Z}\left[\frac{1+i\sqrt{19}}{2}\right]=5\mathbb{Z}\left[\frac{1+i\sqrt{19}}{2}\right] \cdot 7\mathbb{Z}\left[\frac{1+i\sqrt{19}}{2}\right].$$ 
Using Quadratic Reciproicty Law and other properties of Legendre's symbol we obtain that Legendre's symbols $\left(\frac{\Delta_{K}}{7}\right)=\left(\frac{-19}{7}\right)=1$ and $\left(\frac{\Delta_{K}}{7}\right)=\left(\frac{-19}{5}\right)=1.$
Applying the ramification theorem of a prime positive integer in the ring of algebraic integers of a quadratic field, we obtain that each of the ideals $5\mathbb{Z}\left[\frac{1+i\sqrt{19}}{2}\right]$  and $7\mathbb{Z}\left[\frac{1+i\sqrt{19}}{2}\right]$  totally splits in the ring
$\mathbb{Z}\left[\frac{1+i\sqrt{19}}{2}\right].$ Thus, the unique decomposition of $I$ into the product of prime ideals of the ring $\mathbb{Z}\left[\frac{1+i \sqrt{19}}{2}\right]$ is
$$I= \left(\frac{1+i\sqrt{19} }{2}\right)\cdot  \left(\frac{1-i\sqrt{19} }{2}\right)\cdot  \left(\frac{3+i\sqrt{19} }{2}\right)\cdot  \left(\frac{3-i\sqrt{19} }{2}\right)$$
(see \cite{Murty}). The entropy of the ideal $I$ is $H\left(I\right)=\log 4 - \frac{1}{4}\cdot 4\cdot \log 1 =\log 4.$\\
We obtain the following result.
\begin{proposition}
\label{twodoteight}
 Let $K$ be an algebraic number field and let $\mathcal{O}_{K}$ be its  ring of algebraic integers. Let $I$ be an ideal of the ring $\mathcal{O}_{K}.$ Then:
 \begin{equation}\label{sec3_5}
0\leq H\left(I\right)\leq \log \: \omega\left(I\right),
\end{equation}
where $ \omega\left(I\right)$  is the number of distinct prime divisors of the ideal $I.$
\end{proposition}
\begin{proof}
The proof of this Proposition is similar to the proof of Theorem 2  from the article \cite{Minculete1}.             
\end{proof}
We now return to ring $\mathbb{Z}\left[\frac{1+i\sqrt{19}}{2}\right]$ and consider its ideals $I=5\mathbb{Z}\left[\frac{1+i\sqrt{19}}{2}\right]$ and $J=7\mathbb{Z}\left[\frac{1+i\sqrt{19}}{2}\right].$
The unique decompositions of these ideals into products of prime ideals in the ring $\mathbb{Z}\left[\frac{1+i\sqrt{19}}{2}\right]$ are: 
$$I= \left(\frac{1+i\sqrt{19} }{2}\right)\cdot  \left(\frac{1-i\sqrt{19} }{2}\right), \; J= \left(\frac{3+i\sqrt{19} }{2}\right)\cdot  \left(\frac{3-i\sqrt{19} }{2}\right) $$ (see \cite{Murty}). 
So $\omega(I)=\omega(J)=2,$ $\Omega\left(I\right)=\Omega\left(J\right)=2.$ The divergence $D(I||J)=0.$\\
We obtain the following remark.
\begin{remark}
\label{twodotnine}
Let $K$ be an algebriac number field and let $I,J\neq \left(0\right)$ be two ideals of the ring  $\mathcal{O}_{K}$,  uniquely decomposed as $I=P^{e_{1}}_{1}\cdot P^{e_{2}}_{2}\cdot...\cdot P^{e_{g}}_{g}$ and $J=Q^{e^{'}_{1}}_{1}\cdot Q^{e^{'}_{2}}_{2}\cdot...\cdot Q^{e^{'}_{g}}_{g}$,
with $e_{1}, e_{2},...,e_{g}, e^{'}_{1}, e^{'}_{2},...,e^{'}_{g}$ positive integers, $P_{1}, P_{2},...,P_{g}$ distinct prime ideals of the ring $\mathcal{O}_{K}$ and $Q_{1}, Q_{2},...,Q_{g}$ distinct prime ideals of the ring $\mathcal{O}_{K}.$
If $e_{i}=e^{'}_{i},$ for  $i=\overline{1,g}$,  then $D\left(I||J\right)=D\left(J||I\right)=0$.
\end{remark}
\begin{proof}
The proof follows immediately, applying Definition \ref{twodotone.1}.
\end{proof}

\smallskip

\section{The entropy and exponential divisiors}

\noindent

In the paper \cite{Min_Sav}, Minculete and Savin introduced the notion of \textit{exponential divisor} of an ideal of the ring of integers of an algebraic number field.

Keeping the notations from the previous section, an \textit{exponential divisor} of an ideal $I$ of the ring  $\mathcal{O}_{K}$,  $I=P^{e_{1}}_{1}\cdot P^{e_{2}}_{2}\cdot...\cdot P^{e_{g}}_{g}$ has the following form
$d^{\left(e\right)}_{I}=P^{\beta_{1}}_{1}\cdot P^{\beta_{2}}_{2}\cdot...\cdot P^{\beta_{g}}_{g},$
where $\beta_{1},$ $\beta_{2},$...,$\beta_{g}\in\mathbb{N}^{*},$ with $\beta_{i}$ $|$ $e_{i},$ for $\left(\forall\right)$ $i=\overline{1,g}.$

 Let $\mathbb{J}$ be the set of ideals of the ring $\mathcal{O}_{K}.$ Minculete and Savin (in  \cite{Min_Sav}) extended the functions $\tau,$  $\tau^{\left(e\right)}$ to ideals of the ring $\mathcal{O}_{K}$ so:\\
$\tau : \mathbb{J}\rightarrow \mathbb{C},\tau\left(I\right)=$ the number of divisors of the ideal $I,$ respectively\\
$\tau^{\left(e\right)} : \mathbb{J}\rightarrow\mathbb{C}$, $\tau^{\left(e\right)}\left(I\right)=$ the number of exponential divisors of the ideal $I$. 

Also, they gave in \cite{Min_Sav}, the formulas for $\tau\left(I\right)$  and for  $\tau^{\left(e\right)}\left(I\right),$ namely:
\begin{equation}\label{sec3_1}
\tau\left(I\right)=\left(e_{1}+1\right)\cdot\left(e_{2}+1\right)\cdot...\cdot\left(e_{g}+1\right)                       
\end{equation}
and 
\begin{equation}\label{sec3_2}
\tau^{\left(e\right)}\left(I\right)=\tau\left(e_{1}\right)\cdot \tau\left(e_{2}\right)\cdot...\cdot \tau\left(e_{g}\right)    
\end{equation}
(see Proposition 2.1 from \cite{Min_Sav}).

We obtain the following new inequalities about the entropy of a certain type of divisor and about the entropy of an exponential divisor.
\smallskip\\
\begin{proposition}
\label{fourdotone}
Let $\alpha,\beta$ be two natural numbers with $\alpha,\beta\geq 1$, $p$ a prime number and decomposition in prime factors of $n$ given by $n=\prod_{i=1}^rp_i^{a_i}$ with $a_i\in\mathbb{N}^*$ for all $i\in\{1,...,r\}$ such that $\gcd(n,p)=1$. Then\\
i) If $\alpha\geq\beta\geq\Omega(n)e^{-H(n)}$, then we have $H(np^\alpha)\leq H(np^\beta)$;\\
ii) If $\beta\leq\alpha\leq\Omega(n)e^{-H(n)}$, then we have $H(np^\alpha)\geq H(np^\beta)$.
\end{proposition}
\begin{proof}
We study the difference of the entropy of the numbers $np^\alpha$ and $np^\alpha$, with $\gcd(n,p)=1$ and $\alpha\geq\beta\geq 1$. Using equality \eqref{ms_5}, we deduce
\begin{align*}
H(np^\alpha)-H(np^\beta)&=\frac{(\alpha-\beta)(\Omega(n)\log\Omega(n)-\Omega(n)H(n))}{(\Omega(n)+\alpha)(\Omega(n)+\beta)}+\log\frac{\Omega(n)+\alpha}{\Omega(n)+\beta}\\
&+\frac{\beta\log\beta}{\Omega(n)+\beta}-\frac{\alpha\log\alpha}{\Omega(n)+\alpha}\\
&=\frac{(\alpha-\beta)\sum_{i=1}^{r}\alpha_{i}\log \alpha_{i}}{(\Omega(n)+\alpha)(\Omega(n)+\beta)}+\log\frac{\Omega(n)+\alpha}{\Omega(n)+\beta}\\
&+\frac{\beta\log\beta}{\Omega(n)+\beta}-\frac{\alpha\log\alpha}{\Omega(n)+\alpha}.
\end{align*}
We denote $\sum_{i=1}^{r}\alpha_{i}\log \alpha_{i}=A\geq 0$ and $\Omega(n)=t\geq 1$ and we take the function $f:[\beta,\infty)\to\mathbb{R}$ defined by $f(\alpha)=\frac{(\alpha-\beta)A}{(t+\alpha)(t+\beta)}+\log\frac{t+\alpha}{t+\beta}+\frac{\beta\log\beta}{t+\beta}-\frac{\alpha\log\alpha}{t+\alpha}.$ But, $f'(\alpha)=\frac{A-t\log\alpha}{(t+\alpha)^2}=\frac{t}{(t+\alpha)^2}(\frac{A}{t}-\log\alpha)=\frac{\Omega(n)}{(\Omega(n)+\alpha)^2}(\log\Omega(n)-H(n)-\log\alpha),$ which is equivalent to $f'(\alpha)=\frac{\Omega(n)}{(\Omega(n)+\alpha)^2}\log\frac{\Omega(n)e^{-H(n)}}{\alpha}.$
In the first case $\beta\geq\Omega(n)e^{-H(n)}$, then we deduce $f'(\alpha)\leq 0$, so, the function $f$ is decreasing. Therefore for $\alpha\geq\beta$ we have $f(\alpha)\leq f(\beta)=0$, which means that $H(np^\alpha)\leq H(np^\beta)$. The second case $\beta\leq\alpha\leq\Omega(n)e^{-H(n)}$, then $f'(\alpha)\geq 0$, so, the function $f$ is increasing. Since $\alpha\geq\beta$ we deduce $f(\alpha)\leq f(\beta)=0$, which means that $H(np^\alpha)\geq H(np^\beta)$. Consequently the statement is true.
\end{proof}
\begin{remark}
\label{remark}
 Let $n$ be a positive integer,  $n\geq2,$ $n=p^{\alpha_{1}}_{1}p^{\alpha_{2}}_{2}...p^{\alpha_{r}}_{r}$, where $p_{1}, p_{2},...,p_{r}$ are distinct prime positive integers and $\alpha_{1}, \alpha_{2},...,\alpha_{r}\in\mathbb{N}^{*}.$ 
 If $n$ is a square free integer, then $\alpha_{1}= \alpha_{2}=...=\alpha_{r}=1.$ So, $\beta\geq 1=\Omega(n)e^{-H(n)}$, then we have $H(np^\alpha)\leq H(np^\beta)$. Therefore, for any $d_{e}$ be an exponential divisor of $m=np^{\alpha}$, the following inequality hold:
$$H\left(d_{e}\right)\geq H\left(m\right).$$
If $n$ is not a square free integer, then for an exponential divisor $d_{e}$ of $m=np^{\alpha}$ given by $d_{e}=np^{\beta}$, the following inequality hold:
$$H\left(d_{e}\right)\leq H\left(m\right),$$
when $\beta\leq\alpha\leq\Omega(n)e^{-H(n)}$.
\end{remark}

\smallskip

Considering now a positive integer $n$ of the form  $n=p^{\alpha_{1}}_{1}\cdot p^{\alpha_{2}}_{2}\cdot...\cdot p^{\alpha_{r}}_{r}$, where $p_{1}, p_{2},...,p_{r}$ are distinct prime positive integers and $\alpha_{1}, \alpha_{2},...,\alpha_{r}\in$ $\left\{1, 2\right\}$ and an arbitrary exponential divisor $d_{e}$ of $n$, we are interested in the relationship between the entropy $ H\left(n\right)$ and the entropy $H\left(n\right).$

For example, when $r=2$ and $n=12=2^{2}\cdot 3,$ $ H\left(12\right)=\log \: 3 -  \frac{2}{3}\cdot \log \: 2.$ $d_{e}=6$ is an exponential divisor of $12$ and  $ H\left(6\right)=  \log \: 2,$ so $ H\left(12\right)\leq H\left(6\right)$ (we are in the case iii) from Proposition \ref{fourdotone}).

Another example: when $r=3$ and $n=180=2^{2}\cdot 3^{2}\cdot 5,$ $ H\left(180\right)=\log \: 5 -  \frac{4}{5}\cdot \log \: 2.$ $d_{e}=60=2^{2}\cdot 3\cdot 5$ is an exponential divisor of $180$ and  $ H\left(60\right)=  \log \: 4 -  \frac{1}{2}\cdot \log \: 2,$ so $ H\left(60\right)\leq H\left(180\right)$ (we are in the case ii). from Proposition \ref{fourdotone}).

\begin{corollary}
\label{foudottwo}
 Let $n$ be a positive integer,  $n\geq2,$ $n=p^{\alpha_{1}}_{1}\cdot p^{\alpha_{2}}_{2}\cdot...\cdot p^{\alpha_{r}}_{r}$, where $r\geq3,$ $p_{1}, p_{2},...,p_{r}$ are distinct prime positive integers and $\alpha_{1}, \alpha_{2},...,\alpha_{r}\in$ $\left\{1, 2\right\}.$ Then, for any $d_{e}$ be an exponential divisor of $n,$ the following inequality holds:
$$H\left(d_{e}\right)\leq H\left(n\right).$$
\end{corollary}
\begin{proof}
Let $d_{e}$ be an exponential divisor of $n,$ then $d_{e}=p^{\beta_{1}}_{1}\cdot p^{\beta_{2}}_{2}\cdot...\cdot p^{\beta_{r}}_{r}$, where  $\beta_{1}, \beta_{2},$
$...,\beta_{r}\in\mathbb{N}^{*},$  with $ \beta_{i}\in\left\{1, 2\right\},$  $ \beta_{i} |  \alpha_{i}$ $\left(\forall\right)$ $i=\overline{1,r}.$ 
Without diminishing the generality, we assume that there exists $s\in\mathbb{N}^{*},$ $s\leq r$ such that $\alpha_{1}=...=\alpha_{s}=2$ and $\alpha_{s+1}=...=\alpha_{r}=1$ and 
 there exists $l\in\mathbb{N}^{*},$ $l\leq s$ such that $\beta_{1}=...=\beta_{l}=2$ and $\beta_{l+1}=...=\beta_{r}=1.$ Then $n=p^{2}_{1}\cdot p^{2}_{2}\cdot...\cdot p^{2}_{l}\cdot....\cdot p^{2}_{s}\cdot p_{s+1}\cdot...\cdot p_{r}$ and  $d_{e}= p^{2}_{1}\cdot p^{2}_{2}\cdot...\cdot p^{2}_{l}\cdot  p_{l+1}\cdot  ...\cdot p_{s} \cdot  p_{s+1}\cdot...\cdot p_{r}.$
Repeatedly applying Proposition \ref{fourdotone} or Remark \ref{remark}, we have:
$$ H\left(n\right)=H\left(p^{2}_{1}\cdot p^{2}_{2}\cdot...\cdot p^{2}_{l}\cdot....\cdot p^{2}_{s-1}\cdot p^{2}_{s}\cdot p_{s+1}\cdot...\cdot p_{r}\right)\leq $$
$$\leq H\left(p^{2}_{1}\cdot p^{2}_{2}\cdot...\cdot p^{2}_{l}\cdot....\cdot p^{2}_{s-1}\cdot p_{s}\cdot p_{s+1}\cdot...\cdot p_{r}\right)\leq ...\leq  $$
$$\leq H\left(p^{2}_{1}\cdot p^{2}_{2}\cdot...\cdot p^{2}_{l}\cdot  p_{l+1}\cdot...\cdot p_{s}\cdot p_{s+1}\cdot...\cdot p_{r}\right)=H\left(d_{e}\right).$$
\end{proof}

\begin{corollary}
\label{foudotthree}
Let $K$ be an algebraic number field and let $\mathcal{O}_{K}$ be the ring of integers of $K.$  Let $I$ be an ideal of the ring $\mathcal{O}_{K},$  $I=P^{e_{1}}_{1}\cdot P^{e_{2}}_{2}\cdot...\cdot P^{e_{g}}_{g},$ where $g\geq3,$ $P_{1},$ $P_{2},$..., $P_{g}$ are distinct prime ideals of the ring $\mathcal{O}_{K}$ 
and $e_{1}, e_{2},...,e_{r}\in\left\{1, 2\right\}.$  Then, for any  $d^{\left(e\right)}_{I}$ an exponential divisor of the ideal $I,$ the following inequality holds:
$$    H\left(d^{\left(e\right)}_{I}\right) \leq H\left(I\right).$$

\end{corollary}
\begin{proof}
The proof is similar to the proof of Corollary  \ref{foudottwo}.
\end{proof}

We consider the following example of calculating the divergence of two ideals.\\
Let the quadratic field $\mathbb{Q}\left(i\right)$ (where $i^{2}=-1$). Its ring of algebraic integers is  $\mathbb{Z}\left[i\right]$, which is an Euclidean ring with respect to the norm $N$ ($N: \mathbb{Z}\left[i\right]\rightarrow \mathbb{N}$, $N\left(a+b\cdot i\right)=a^{2}+b^{2}$), so it is a principal ring, so it is a Dedeking ring.
We consider the ideal $I=90\mathbb{Z}\left[i\right].$ Since $N\left(1+ i\right)=2,$ $N\left(1+2 i\right)=N\left(1-2 i\right)=5$ and $2$ and $5$ are prime positive integers, it results that $\left(1+ i\right)\mathbb{Z}\left[i\right],$  $\left(1+2 i\right)\mathbb{Z}\left[i\right],$ $\left(1-2 i\right)\mathbb{Z}\left[i\right]$ are prime ideals. Taking into account that $2$ and $\left(1+ i\right)^{2}$ are associated in divisibility in the ring  $\mathbb{Z}\left[i\right]$, so the ideals $2\mathbb{Z}\left[i\right]$ and  $\left(1+ i\right)^{2}\cdot\mathbb{Z}\left[i\right]$ are equal and the fact that any prime positive integer congruent to $3$ mod $4$ is inert in the ring 
 $\mathbb{Z}\left[i\right]$ (see \cite{ireland}), we obtain that the unique decomposition of ideal  $I=90\mathbb{Z}\left[i\right]$ into the product of prime ideals of the ring $\mathbb{Z}\left[i\right]$ is:

$$I=90\mathbb{Z}\left[i\right]=\left(1+ i\right)^{2}\cdot \mathbb{Z}\left[i\right]\cdot 3^{2}\cdot \mathbb{Z}\left[i\right]\cdot \left(1+2 i\right) \mathbb{Z}\left[i\right]\cdot \left(1+2 i\right) \mathbb{Z}\left[i\right].$$
Let $\gamma(I)$ be the lowest exponential divisor of ideal $I$, that is 
$$\gamma(I)=\left(1+ i\right)\mathbb{Z}\left[i\right]\cdot 3\mathbb{Z}\left[i\right]\cdot \left(1+2 i\right) \mathbb{Z}\left[i\right]\cdot \left(1+2 i\right) \mathbb{Z}\left[i\right].$$
$\omega(I)=\omega(\gamma(I))=4,$ $\Omega\left(I\right)=6,$ $\Omega\left(\gamma(I)\right)=4.$ We have:

$$ D(I||\gamma(I))=\log\frac{2}{3}+  \frac{2}{3}\cdot \log 2; \; H\left(I\right)= \log 6 - \frac{2}{3}\cdot \log 2.  $$
So, $ D(I||\gamma(I))+ H\left(I\right)=\log4=\log\omega(\gamma(I))=H\left(\gamma(I)\right).$\\
\smallskip\\
We now generalize the first part of Remark  \ref{remark2.2},  for ideals.
\begin{remark}
\label{remark4.4.}
Let $K$ be an algebraic number field and let $I\neq \left(0\right)$ be an ideal of the ring  $\mathcal{O}_{K}$, uniquely decomposed as $I=P^{e_{1}}_{1}\cdot P^{e_{2}}_{2}\cdot...\cdot P^{e_{g}}_{g},$ with $e_{1}, e_{2},...,e_{g}$ positive integers and $P_{1}, P_{2},...,P_{g}$ are distinct prime ideals of the ring $\mathcal{O}_{K}$.
Let $\gamma(I)=P_{1}\cdot P_{2}\cdot...\cdot P_{g}$ be the lowest exponential divisor of ideal $I$.
Then, the following equality holds: 
$$ D(I||\gamma(I))=H(\gamma(I))- H\left(I\right)=\log\omega(I)- H\left(I\right).
$$
\end{remark}
The proof of this remark is similar to the proof of Remark  \ref{remark2.2}.

\section{Conclusions} 
In the paper \cite{MS} we introduced the notion of entropy of an ideal generalizing the notion of entropy of a positive integer, introduced by Minculete and Pozna in \cite{Minculete1} and we found some properties of this entropy. Also, we found some properties of the entropy of the ideals of the form 
$p\mathcal{O}_{K},$ where $\mathcal{O}_{K}$ is the ring of integers of an algebraic number field $K$ and $p$ is a prime positve integer. In this paper we generalize a certain inequality satisfied by the entropy of ideals of the type $p\mathcal{O}_{K}$ to the entropy of arbitrary ideals $I$ of the ring $\mathcal{O}_{K}.$ We also introduce the notion of divergence of two arbitrary ideals of the ring  $\mathcal{O}_{K}$  and we find some properties of it, generalizing some properties involving the divergence of two positive integers. We notice that many of the inequalities satisfied by the entropy of a positive integer or by the divergence of two positive integers are carried over in the same way to ideals, but we also obtained other interesting results for the entropy of an ideal of the type $p\mathcal{O}_{K}$ (with  $p$ a prime positve integer ), respectively about the divergence of two ideals (see the upper bound of the inequality in Proposition \ref{twodotfive} and  see Proposition \ref{twodotseven}).

In addition, we find new inequalities satisfied by the entropy $H$ of a positive integer and other results involving both the notion of entropy and the notion of divergence.

In the last section we find some inequalities, involving the  entropy  $H$ of an  exponential divisor of a positive integer, respectively the  entropy  $H$ of an  exponential divisor of an ideal.

In the future, we will look for other connections of entropy within ideals, studying a possible generalization of existing entropy types for natural numbers or for ideals.

\end{document}